\documentclass[a4paper,12pt]{article}
\usepackage{amsfonts,amsmath}
\usepackage{amssymb}
\usepackage[mathscr]{eucal}
\usepackage{amsthm}
\usepackage{amscd}
\usepackage[T2A]{fontenc}
\usepackage{graphics}
\usepackage[dvips]{color}
\usepackage[cp1251]{inputenc}

\newcounter{rem}[section]

\def\text{\mbox}

\def\varkappa{\kappa}
\def\suml{\sum\limits}
\def\maxl{\max\limits}
\def\minl{\min\limits}
\def\intl{\int\limits}
\def\supl{\sup\limits}
\def\infl{\inf\limits}

\def\prodl{\prod\limits}

\def\dpfrac{\displaystyle\frac}

\newcommand{\df}{\stackrel{\textrm{def}}{=}}

\def\text{\mbox}

\def\dzeta{\zeta}
\def\suml{\sum\limits}
\def\maxl{\max\limits}
\def\minl{\min\limits}
\def\intl{\int\limits}
\def\supl{\sup\limits}
\def\infl{\inf\limits}
\def\liminfl{\liminf\limits}

\def\prodl{\prod\limits}
\def\dip{\displaystyle}

\def\dpfrac{\dip\frac}
\title{Asymptotically optimal detection of changes in stochastic models with switching regimes}
\author{Brodsky B.E., Darkhovsky B.S.}
\date{}

\begin{document}
\sloppy \maketitle

\begin{abstract}
This paper deals with the problem of asymptotically optimal
detection of changes in regime-switching stochastic models. We
need to divide the whole obtained sample of data into several
sub-samples with observations belonging to different states of a
stochastic models with switching regimes. For this purpose, the
idea of reduction to a corresponding change-point detection
problem is used. Both univariate and multivariate switching models
are considered. For the univariate case, we begin with the study
of binary mixtures of probabilistic distributions. In theorems 1
and 2 we prove that type 1 and type 2 errors of the proposed
method converge to zero exponentially as the sample size tends to
infinity. In theorem 3 we prove that the proposed method is
asymptotically optimal by the rate of this convergence in the
sense that the lower bound in the a priori informational
inequality is attained for our method. Several generalizations to
the case of multiple univariate mixtures of probabilistic
distributions are considered. For the multivariate case, we first
study the general problem of classification of the whole  array of
data into several sub-arrays of observations from different
regimes of a multivariate stochastic model with switching states.
Then we consider regression models with abnormal observations and
switching sets of regression coefficients. Results of a detailed
Monte Carlo study of the proposed method for different stochastic
models with switching regimes are presented.

\end{abstract}

\bigskip
{\large \bf 1. Introduction}

In this paper the problem of the retrospective detection of
changes in stochastic models with switching regimes is considered.
Our main goal is to propose asymptotically optimal methods for
detection and estimation of possible 'switches', i.e. random and
transitory departures from prevailing stationary regimes of
observed stochastic models.

First, let us mention previous important steps into this field.
Models with switching regimes have a long pre-history in
statistics (see, e.g., Lindgren (1978)). A simple switching model
with two regimes has the following form:
$$\begin{array}{ll}
& Y_t=X_t\beta_1+u_{1t} \quad \text{for the 1st regime  } \\
& Y_t=X_t \beta_2+u_{2t} \quad \text{for the 2nd regime }.
\end{array}
$$

For models with endogenous switchings usual estimation techniques
for regressions are not applicable. Goldfeld and Quandt (1973)
proposed {\it regression models with Markov switchings}. In these
models probabilities of sequential switchings are supposed to be
constant. Usually they are described by the matrix of
probabilities of switchings between different states.

Another modification of the regression models with Markov
switchings was proposed by Lee, Porter (1984). The following
transition matrix was studied:
$$
\Lambda=[p_{ij}]_{i,j=0,1}, \quad p_{ij}=P\{I_t=j|I_{t-1}=i\}.
$$

Lee and Porter (1984) consider an example with railway transport
in the US in 1880-1886s which were influenced by the cartel
agreement. The following regression model was considered:
$$
log P_t=\beta_0+\beta_1 X_t+\beta_2 I_t+u_t,
$$
where $I_t=0$ or $I_t=1$ in dependence of 'price wars' in the
concrete period.

Cosslett and Lee (1985) generalized the model of  Lee and Porter
to the case of serially correlated errors $u_t$.

Many economic time series occasionally exhibit dramatic breaks in
their behavior, assocoated with with events such as financial
crises (Jeanne and Mason, 2000; Cerra, 2005; Hamilton, 2005) or
abrupt changes in government policy (Hamilton, 1988; Sims and Zha,
2004; Davig, 2004). Abrupt changes are also a prevalent feature of
financial data and empirics of asset prices (Ang and Bekaert,
2003; Garcia, Luger, and Renault, 2003; Dai, Singleton, and Wei,
2003).

The functional form of the 'hidden Markov model' with switching
states can be written as follows:
$$
y_t=c_{s_t}+\phi y_{t-1}+\epsilon_t, \eqno (i)
$$
where $s_t$ is a random variable which takes the values $s_t=1$
and $s_t=2$ obeying a two-state Markov chain law:
$$
Pr(s_t=j|s_{t-1}=i,s_{t-2}=k,\dots,y_{t-1},y_{t-2},\dots)=Pr(s_t=j|s_{t-1}=i)=p_{ij}.
\eqno(ii)
$$

A model of the form (1-2) with no autoregressive elements
($\phi=0$) appears to have been first analyzed by Lindgren (1978)
and Baum, et al. (1980). Specifications that incorporate
autoregressive elements date back in the speech recognition
literature to Poritz (1982), Juang and Rabiner (1985), and Rabiner
(1989). Markov-switching regressions were first introduced in
econometrics by Goldfeld and Quandt (1973), the likelihood
function for which was first calculated by Cosslett and Lee
(1985). General characterizations of moment and stationarity
conditions for Markov-switching processes can be found in
Tjostheim (1986), Yang (2000), Timmermann (2000), and Francq and
Zakoian (2001).

A useful review of modern approaches
 to estimation in Markov-switching models can be found in Hamilton (2005).

However, the mechanism of Markov chain modeling is far not unique
in statistical description of dependent observations. Besides
Markov models, we can mention martingale and copula approaches to
dealing with dependent data, as well as description of statistical
dependence via different coefficients of 'mixing'. All of these
approaches are interrelated and we must choose the most
appropriate method for the concrete problem. In this paper we
choose the mixing paradigm for description of statistical
dependence.

Remark that $\psi$-mixing condition is imposed below in this paper
in order to obtain the exponential rate of convergence to zero for
type 1 and type 2 error probabilities (see theorems 1 and 2
below). Another alternative was to assume $\alpha$-mixing property
which is always satisfied for aperiodic and irreducible
countable-state Markov chains (see Bradley (2005)). Then we can
obtain the hyperbolic rate of convergence to zero for type 1 and
type 2 error probabilities. For the majority of practical
applications, it is enough to assume $r$-dependence (for a certain
finite number of lags $r\ge 1$) of observations and state
variables. Then all proofs become much shorter.

Now let us mention some important problems which lead to
stochastic models with switching regimes.

\pagebreak {\it Splitting mixtures of probabilistic distributions}

In the simplest case we suppose that the d.f. of observations has
the following form:
$$
F(x)=(1-\epsilon)F_0(x)+\epsilon F_1(x),
$$
where $F_0(x)$ is the d.f. of ordinary observations; $F_1(x)$ is
the d.f. of abnormal observations; $0\le \epsilon <1$ is the
probability of obtaining an abnormal observation.

We need to test the hypothesis of statistical homogeneity (no
abnormal observations) of an obtained sample
$X^N=\{x_1,x_2,\dots,x_N\}$. If this hypothesis is rejected then
we need to estimate the share of abnormal observations
($\epsilon$) in the sample and to classify this sample into
sub-samples of ordinary and abnormal observations.

\bigskip
{\it Estimation for regression models with abnormal observations}

The natural generalization of the previous model is the regression
model with abnormal observations
$$
Y=X\beta+\epsilon,
$$
where $Y$ is the $n\times 1$ vector of dependent observations; $X$
is the $n\times k$ matrix of predictors; $\beta$ is $k\times 1$
vector of regression coefficients; $\epsilon$ id the $n\times 1$
vector of random noises with the d.f. of the following type:
$$
f_{\epsilon}(x)=(1-\delta) f_0(x)+\delta f_1(x),
$$
where $0\le \delta <1$ is the probability to obtain an abnormal
observation; $f_0(x)$ is the density function of ordinary
observations; $f_1(x)$ is the density function of abnormal
observations. For example, in the model with Huber's contamination
[Huber, 1985]: $f_0(\cdot)={\cal N}(0,\sigma^2),\;f_1(\cdot)={\cal
N}(0,\Lambda^2)$.

\bigskip
{\it Estimation for regression models with changing coefficients}

Regression models with changing coefficients is another
generalization of the contamination model. Suppose a baseline
model is described by the following regression:
$$
Y=X\alpha+\xi,
$$
where the mechanism of a change is purely random:
$$
\alpha=\left \{
\begin{array}{ll}
& \beta \text{ with the probability   }   1-\epsilon \\
& \gamma \text{ with the probability   }   \epsilon,
\end{array}
\right.
$$
and $\beta\ne \gamma$.

We need again to test the hypothesis of statistical homogeneity of
an obtained sample and to divide this sample into sub-samples of
ordinary and abnormal observations if the homogeneity hypothesis
is rejected.

The goal of this paper is to propose methods which can solve these
problems effectively. Theoretically, we mean estimation of type 1
and type 2 errors in testing the statistical homogeneity
hypothesis and with estimation of contaminations parameters in the
case of rejectiong this hypothesis. Practically, we propose
procedures for implementation of these methods for univariate and
multivariate models.

The structure of this paper is as follows. First, we consider
univariate models with switching effects. For binary mixtures of
probabilistic distributions we prove theorem 1 about exponential
convergence to zero of type 1 error in classification (to detect
switches for a statistically homogenous sample) as the sample size
$N$ tends to infinity; theorem 2 about exponential convergence to
zero of type 2 error (vice versa, to accept stationarity
hypothesis for a sample with switches); and theorem 3 which
establishes the lower bound for the error of classification for
binary mixtures. From theorems 2 and 3 we conclude that the
proposed method is asymptotically optimal by the order of
convergence to zero of the classification error.

Different generalizations of the proposed method for the case of
univariate models with multiple switching regimes and for
multivariate models with switching regimes are considered. Results
of a detailed Monte Carlo study of the proposed method for
different stochastic models with switching regimes are presented.

\vspace{1.5cm}

 {\large \bf 2. Univariate models}

\bigskip
{\large \bf 2.1. Binary mixtures}

\bigskip
 {\bf 2.1.1. Problem statement and description of the detection/estimation method}

Suppose the d.f. of the observations is the binary mixture
$$
f(x)=(1-\epsilon)f_0(x)+\epsilon f_0(x-h),
$$
where $\epsilon, h$ are unknown.

The problem is to estimate parameters $\epsilon,h$ by the sample $X^N=\{x_i\}_{i=1}^N$, where all $x_i$ has the same d. f. $f(\cdot)$.

An ad hoc method of estimation of these parameters is as follows: ordinary and 'abnormal' observations are heuristically classified to two
sub-samples and the estimate $\hat\epsilon$ is computed as the share of the size of sub-sample of abnormal observations in the whole sample
size. Clear, this method is correct only for large values of $h$. However, this idea of two sub-samples can be used in construction of more
subtle methods of estimation.

The estimation method is as follows:

1) From the initial sample $X^N$  compute the estimate of the mean value:
$$
\theta_N=\dpfrac 1N \suml_{i=1}^n\,x_i
$$

2) Fix the parameter $b>0$ and classify observations as follows: if an observation falls into the interval $(\theta_N-b,\theta_N+b)$, then we
place it into the sub-sample of ordinary observations, otherwise - to the sub-sample of abnormal observations.

3) Then for each $b>0$ we obtain the following decomposition of the sample $X^N$ into two sub-samples
$$\begin{array}{ll}
& X_1(b)=\{\tilde x_1,\tilde x_2,\dots,\tilde x_{N_1}\},\quad |\tilde x_i-\theta_N|< b, \\
& X_2(b)=\{\hat x_1,\hat x_2,\dots,\hat x_{N_2}\},\quad |\hat x_i-\theta_N|\ge b
\end{array}
$$
Denote by $N_1=N_1(b),\,N_2=N_2(b),\,N=N_1+N_2$ the sizes of the sub-samples $X_1$ and $X_2$, respectively.

The parameter  $b$  is chosen so that the sub-samples $X_1$ and $X_2$ are separated in the best way. For this purpose, consider the following
statistic:
$$
\Psi_N(b)=\dpfrac 1{N^2}(N_2 \suml_{i=1}^{N_1}\,\tilde x_i-N_1\suml_{i=1}^{N_2}\,\hat x_i).
$$

4) Define the boundary $C>0$ and compare it with the value $J=\max |\Psi_N(b)|$ on the set $b>0$. If $J\le C$ then we accept the hypothesis
$H_0$ about the absence of abnormal observations; if, however, $J>C$ then the hypothesis $H_0$ is rejected and the estimates of the parameters
$\epsilon$ and $h$ are constructed. Remark that our primary goal is to separate ordinary and abnormal observations in the sample. Evidently,
classification errors must be small and therefore we have to require some kind of convergence of the estimate $\hat\epsilon_N$ to its true value
$\epsilon$.

5) If $J>C$ then define the number $b_N^*$:
$$
b_N^* \in \arg\max_{b>0}|\Psi_N(b)|.
$$

Then
$$
\epsilon_N^*=N_2(b_N^*)/N,\quad h_N^*=\theta_N / \epsilon_N^*.
$$
are the nonparametric estimates for $\epsilon$ and $h$, respectively.

\bigskip

In the general case for construction of unbiased and consistent estimates of $\epsilon$ and $h$ we can use the following relationships:
$$\begin{array}{ll}
& \hat\epsilon_N \,\hat h_N=\theta_N \\
& \dpfrac {1-\hat\epsilon_N}{\hat\epsilon_N}=\dpfrac {f_0(\theta_N-b_N^*-\hat h_N)-f_0(\theta_N+b_N^*-\hat
h_N)}{f_0(\theta_N+b_N^*)-f_0(\theta_N-b_N^*)}.
\end{array}
$$

We will show that, under some conditions, the estimates $\hat\epsilon_N$ and $\hat h_N$ tend almost surely to the true values $\epsilon$ and $h$
as $N\to\infty$. The sub-sample of abnormal observations is  $X_2(b_N^*)$.

\vspace{0.5 cm} {\bf 2.1.2. Main results}

Let us formulate the main assumptions.

A1.\textbf{Mixing conditions}

On the probability space $(\Omega,\mathfrak {F},\mathbf{P})$
let $\mathcal{H}_{1}$ and
$\mathcal{H}_{2}$ be two $\sigma $-algebras from $\mathfrak{F}$. Consider the following measure of dependence between $\mathcal{H}_{1}$ and
$\mathcal{H}_{2}$:
$$
 \psi (\mathcal{H}_{1}, \mathcal{H}_{2})  = \sup_{A \in
\mathcal{H}_{1},B \in \mathcal{H}_{2}, \mathbf{P}(A) \mathbf{P}(B)\neq 0} \Big\vert \frac{\mathbf{P}(AB)}{\mathbf{P}(A) \mathbf{P}(B)}-1
\Big\vert
$$

Suppose $\{y_n\},\,n\ge 1$ is a sequence of random variables defined on ($\Omega,\mathfrak{F},\mathbf{P}$). Denote by $\mathfrak{F}^{t}_{s}=\sigma
\{y_{i}: s\le i\le t\}, 1\le s\le t< \infty$ the minimal $\sigma $-algebra generated by random variables $y_{i}, s \le i \le t$. Define
$$
 \psi (n)  = \sup_{t\ge 1} \psi (\mathfrak{F}^{t}_{1},
\mathfrak{F}^{\infty }_{t+n})
$$

We say that random sequence $\{y_n\}$ satisfies the $\psi $-\emph{mixing condition} if the function $\psi(n)$ (which is also called the
 \emph{$\psi $-mixing coefficient}) tends to zero as $n$ goes to infinity.

The $\psi$-mixing condition is satisfied in most practical cases. In particular, for a Markov chain (not necessarily stationary), if $\psi(n)<1$
for a certain $n$, then $\psi(k)$ goes to zero at least exponentially as $k\to\infty$ (see Bradley, 2005, theorem 3.3).

A2.\textbf{Cramer condition}

We say that the sequence $\{y_n\}$ satisfies the \emph{uniform Cramer condition} if there exists $T>0$ such that for each $i$, $\mathbf{E}\exp(ty_i) < \infty$ for $|t|<T$.

For a centered sequence $\{y_n\}$ this condition is equivalent to
the following (see Petrov, 1987): there exist $g>0,\,H>0$ such
that
$$
\mathbf{E}e^{ty_n}\le e^{\frac 12 g t^2} , \qquad |t|\le H,
$$
for all $n=1,2,\dots$.

\emph{We assume that conditions A1 and A2 hold true everywhere in the paper.}

For any $x>0$ let us choose the number $\gamma(x)$ from the following condition:
$$
\ln(1+\gamma(x))=\left \{
\begin{array}{ll}
& \dpfrac {x^2}{4g}, \qquad x\le gH \\
& \dpfrac {xH}4, \qquad x>gH,
\end{array}
\right.
$$
where $g,H$ are taken from the uniform Cramer condition.

For the chosen $\gamma(x)$, let us find such integer $\phi_0(x)\ge
1$ from the $\psi$-mixing condition that $\psi(l)\le\gamma(x)$ for
$l\ge \phi_0(x)$.

Below we denote by $\mathbf{P}_0(\mathbf{E}_0),\,\,\mathbf{P}_{\epsilon}(\mathbf{E}_{\epsilon})$ measure (mathematical expectation) of the sequence $X^N$ under the condition $\epsilon=0$ or $h=0$ (no 'abnormal' observations) and under the condition $\epsilon h\not=0$.

In the following theorem the exponential upper estimate for type 1 error  is obtained for the proposed method.

\bigskip
{\bf Theorem 1.}

Let $\epsilon=0$. Suppose the d.f. $f_0(\cdot)$ is symmetric w.r.t. zero. Then for any $C>0$ the following estimate holds:
$$
\mathbf{P}_0\{\maxl_{b>0}\,|\Psi_N(b)| >C\} \le 4\phi_0(CN/2) \exp(-L(C)N),
$$
where $L(C)=min\left(\dpfrac {HC}{8\phi_0(CN/2)},\;\dpfrac {C^2}{16\phi_0^2(CN/2)g}\right)$, the constants $g,H$ are taken from the uniform Cramer  condition.

\bigskip
The proof of theorem 1 is given in the Appendix.

\bigskip
Now consider characteristics of this method in case $\epsilon h\ne 0$. Here we again assume that $\mathbf{E}_0\,x_i=0, \,i=1,\dots,N$.

Put (for any fixed $\epsilon, h$)
$$
\begin{array}{ll}
&r(b)=\intl_{\epsilon h-b}^{\epsilon h+b}\,f(x)x dx,\quad d(b)=\intl_{\epsilon h-b}^{\epsilon h+b}\,f(x) dx\\[2mm]
&\Psi(b)=r(b)-\epsilon h d(b)
\end{array}
$$
and consider the equation
$$
f(\epsilon h+b)=f(\epsilon h-b). \eqno(1)
$$

In the following theorem type 2 error is studied.

\vspace{1cm} {\bf Theorem 2.}

Suppose all assumptions of theorem 1 are satisfied and there exists $r^*=\supl_b\,r(b)$. Suppose also that $f^{''}(\cdot)\not=0$ and continuous. Then for $0< C< \maxl_b\,|\Psi(b)|$
we have

1)
$$
\mathbf{P}_{\epsilon}\{\maxl_b\,|\Psi_N(b)|\le C\}\le 4\phi_0(CN/2+r^*)\,\exp(-L(\delta) N)
$$
where $\delta=\maxl_b\,|\Psi(b)|-C >0,\;L(\delta)=min(\dpfrac
{\delta^2}{16\phi_0^2g},\dpfrac {H\delta}{8 \phi_0})$.

2) Suppose, moreover, that equation (1) has a unique root $b^{*}$ (for any fixed $\epsilon, h$). Then

 $b_N^{*}\to b^{*}$ $\mathbf{P}_{\epsilon}$-a.s. as $N\to\infty$;

3) The estimates $\hat\epsilon_N, \,\hat h_N$ converge $\mathbf{P}_{\epsilon}$-a.s. to the true values of the parameters $\epsilon, h$, respectively, as $N\to\infty$.

\bigskip
The proof of theorem 2 is given in the Appendix.

\bigskip
{\bf 2.1.3. Recommendations for the choice of the threshold $C$}

For practical applications of the above obtained results we need to know the threshold $C$.

a) In order to compute this threshold, at least one training sample without switchings is needed.

b) For this sample we compute the threshold  $C$ from the following empirical formula which follows from theorem 1:
$$
C=C(N)\sim \sigma \sqrt{\dpfrac {\phi_0(\cdot)\,|\ln\alpha|}{N}},
$$
where $N$ is the sample size, $\sigma^2$ is the variation of $\phi_0$-dependent observations and $\alpha$ is the 1st type error level.

In other words we compute the dispersion $\sigma^2$ of observations and the integer $\phi_0$ (by the first zero lag of the autocorrelation
function of the training sample). Then we compute the threshold $C$.

Let us give one example which explains how to do it in practice.

Consider the following model (without switchings)
$$
x(n)=\rho x(n-1)+\sigma\,\xi_n, \qquad n=1,\dots,N,
$$
where $\xi_n$ are i.i.d.r.v.'s with the d.f. $N(0,1)$, and replacing $\phi_0(\cdot)$ by $(1-\rho)^{-1}$.

As a result, the following regression relationship
for the threshold $C$ was obtained:
$$
log(C)=-0.9490-0.4729*log(N)+1.0627*log(\sigma)-0.6502*log(1-\rho)-0.2545*log(1-\alpha).
\eqno(2)
$$

Remark that $R^2=0.978$ for this relationship and its residuals
are stationary at the error level $5 \%$. The elasticity
coefficient for the factor $N$ is close to its theoretical value
$1/2$. The calibration coefficient $\exp(-0.949)=0.3871$ here
depends on the Gaussian d.f. of observations.

We have to note that in practice, we need
to calibrate the above formula for the threshold $C$ using several
homogenous samples.

\bigskip

{\bf Examples}

 The proposed method was tested in the following experiments.

In the first series of tests the following mixture model was studied:
$$
f_{\epsilon}(x)=(1-\epsilon)f_0(x)+\epsilon f_0(x-h), \quad f_0(\cdot)={\cal N}(0,1), \quad 0\le \epsilon <1/2.
$$

First, the critical thresholds of the decision statistic $\max_{b>0} |
\Psi_N(b)|$ were computed. For this purpose we use the
above formula for the threshold $C$ for the values $\alpha=0.95,\;
\rho=0,\;\sigma=1$.The threshold values $C$ in each experiment are
presented in table 1.

\bigskip
{\bf Table 1.}

\bigskip
\hspace{-0.9cm }
\begin{tabular}{|c|c|c|c|c|c|c|c|c|c|c|}
\hline
$N$ & 50 & 100 &  300 & 500 & 800 & 1000 & 1200 & 1500 & 2000 \\
\hline
$\alpha=0.95$ & 0.1681 & 0.1213 & 0.0710 & 0.0534 & 0.044 & 0.0380 & 0.037 & 0.034 & 0.029 \\
\hline
$\alpha=0.99$ & 0.1833 & 0.1410 & 0.0869 & 0.0666 & 0.050 & 0.0471 & 0.0390 & 0.038 & 0.035 \\
\hline
\end{tabular}

\bigskip
In the second series of tests the threshold value for $\alpha=0.95$ was chosen as the critical
threshold $C$ in experiments with non-homogenous samples (for
$\epsilon\ne 0$). For different sample sizes in 5000 independent
trials of each test, the estimate of type 2 error $w_2$ ( i.e. the
frequency of the event $\maxl_b | \Psi_N(b)| <C$ for $\epsilon
>0$) and the estimate $\hat\epsilon$ of the parameter $\epsilon$ were computed. The results are presented in table 2.

\bigskip
{\bf Table 2.}

\bigskip
\begin{tabular}{|c|c|c|c|c|c|c|c|c|}
\hline
$\epsilon=0.1 $ & \multicolumn{4}{|c|} {h=2.0} & \multicolumn{4}{|c|} {h=1.5} \\
\hline
$N$ & 300 & 500 &  800 & 1000 & 800 & 1200 & 2000 & 3000 \\
\hline
$C$ & 0.0710 & 0.0534 & 0.044 & 0.038 & 0.044 & 0.037 & 0.029 & 0.022 \\
\hline
$w_2$ & 0.26 & 0.15 & 0.05 & 0.02 & 0.62 & 0.42 & 0.16 & 0.03 \\
\hline
$\hat\epsilon$ & 0.104 & 0.101 & 0.097 & 0.099 & 0.106 & 0.103 & 0.102 & 0.0985  \\
\hline
\end{tabular}

\vspace{1.5cm}
 {\bf 2.1.4. Asymptotic optimality}

\bigskip
Now consider the question about the asymptotic optimality of the proposed method in the class of all estimates of the parameter $\epsilon$. The
a priori theoretical lower bound for the estimation error of the parameter $\epsilon$ in the model with i.i.d. observations with d.f. $f_{\epsilon}(x)=(1-\epsilon)f_0(x)+\epsilon
f_1(x)$ is given in the following theorem.

\bigskip
{\bf Theorem 3.} Let ${\cal M}_N$ be the class of all estimates of
the parameter $\epsilon$. Then for any $0< \delta < \epsilon$,
$$
\liminfl_{N\to\infty}\infl_{\hat\epsilon_N\in {\cal M}_N}\,\supl_{0< \epsilon <1/2}\,\dpfrac 1N \ln \mathbf{P}_{\epsilon}\{|\hat\epsilon_N-\epsilon|>
\delta \} \ge -\delta^2\,J(\epsilon),
$$
where $J(\epsilon)=\int\,[(f_0(x)-f_1(x))^2/ f_{\epsilon}(x)]\,dx$ is the generalized $\varkappa^2$ distance between densities $f_0(x)$ and
$f_1(x)$ and $\mathbf{P}_{\epsilon}$ is the measure corresponding to the density $f_{\epsilon}(x)$.

{\bf Proof.}

Remark that it suffices to consider consistent estimates of the
parameter $\epsilon$ (for non-consistent estimates the limit in
the left hand of the above inequality is equal to zero). This
class is not empty because of the method proposed in the paper.

Suppose $\hat\epsilon_N$ is any consistent estimate of $\epsilon$
and $0< \delta < \delta^{'}$. Consider the random variable
$\lambda_N=\lambda_N(x_1,\dots,x_N)=\mathbb{I}\{|\hat\epsilon_N-\epsilon|>
\delta \}$, where $\mathbb{I}(A)$ is the indicator of the set $A$.

Then for any $d>0$:
$$
\mathbf{P}_{\epsilon}\{|\hat\epsilon_N-\epsilon| > \delta
\}=E_{\epsilon}\lambda_N\ge \mathbf{E}_{\epsilon} (\lambda_N \mathbb{I}\{f(X^N,
\epsilon+\delta^{'})/ f(X^N,\epsilon) < e^d\}),
$$
where $f(X^N,\epsilon)$ is the likelihood function of the sample
$X^N$ of observations with the density function $f_{\epsilon}(x)$,
i.e.
$$
f(X^N,\epsilon)=\prodl_{i=1}^N\,[(1-\epsilon)f_0(x_i)+\epsilon
f_1(x_i)].
$$

Further,
$$\begin{array}{ll}
& \mathbf{E}_{\epsilon}(\lambda_N\mathbb{I}\{\dpfrac {f(X^N,\epsilon+\delta^{'})}{f(X^N,\epsilon)}<e^d\})\ge \\
& \ge e^{-d}\mathbf{E}_{\epsilon+\delta^{'}}(\lambda_N \mathbb{I}\{f(X^N,\epsilon+\delta^{'})/ f(X^N,\epsilon) < e^d \} \ge \\
& \ge e^{-d}\,(\mathbf{P}_{\epsilon+\delta^{'}}\{|\hat\epsilon_N-\epsilon|>
\delta\}-\mathbf{P}_{\epsilon+\delta^{'}}\{f(X^N,\epsilon+\delta^{'})/
f(X^N,\epsilon)
> e^d\}).
\end{array}
$$

Since $\hat\epsilon_N$ is a consistent estimate,
$\mathbf{P}_{\epsilon+\delta^{'}}\{|\hat\epsilon_N-\epsilon| > \delta \}\to
1$ as $N\to\infty$.

Let us consider the probability
$\mathbf{P}_{\epsilon+\delta^{'}}\{f(X^N,\epsilon+\delta^{'})/
f(X^N,\epsilon)
> e^d\}$. We have
$$\begin{array}{ll}
& \ln\dpfrac {f(X^N,\epsilon+\delta^{'})}{f(X^N,\epsilon)}=\suml_{i=1}^N\,\ln (1+\delta^{'}\dpfrac {f_1(x_i)-f_0(x_i)}{f_{\epsilon}(x_i)})= \\
& =\delta^{'}\suml_{i=1}^N\,\dpfrac
{f_1(x_i)-f_0(x_i)}{f_{\epsilon}(x_i)}+o(\delta^{'}).
\end{array}
$$
On the other hand,
$$
\mathbf{E}_{\epsilon+\delta^{'}}\dpfrac
{f_1(x_i)-f_0(x_i)}{f_{\epsilon}(x_i)}=\delta^{'}\int\,\dpfrac
{(f_1(x_i)-f_0(x_i))^2}{f_{\epsilon}(x_i)}dx_i=\delta ^{'}
J(\epsilon).
$$
Therefore, choosing
$d=N((\delta^{'})^2+\varkappa)J(\epsilon),\;\varkappa=o((\delta^{'})^2)$,
we obtain
$$
\mathbf{P}_{\epsilon+\delta^{'}}\{f(X^N,\epsilon+\delta^{'})/
f(X^N,\epsilon)
> e^d\}\to 0 \text {  as  } N\to\infty.
$$
Thus,
$$
\mathbf{P}_{\epsilon}\{|\hat\epsilon_N-\epsilon| > \delta \} \ge
(1-o(1))\,e^{-N\delta^2\,J(\epsilon)},
$$
or
$$
\liminfl_{N\to\infty}\infl_{\hat\epsilon_N\in {\cal
M}_N}\,\supl_{0< \epsilon <1/2}\,\dpfrac 1N \ln
\mathbf{P}_{\epsilon}\{|\hat\epsilon_N-\epsilon|> \delta \} \ge
-\delta^2\,J(\epsilon),
$$

Theorem 3 is proved.

Comparing results of theorems 2 and 3 we conclude that the proposed method is asymptotically optimal by the order of convergence of the
estimates of a mixture parameters to their true values.

\vspace{1cm} {\bf 2.1.5. Generalizations: non-symmetric distribution functions}

Results obtained in theorems 1 and 2 can be generalized to the case of non-symmetric distribution functions. Suppose the d.f. $f_0(\cdot)$ is
asymmetric w.r.t. zero. Then we can modify the proposed method as follows.

1. From the initial sample $X^N=\{x_1,\dots,x_N\}$ compute the mean value $\theta_N=\dpfrac 1N\,\suml_{i=1}^N\,x_i$ and the sample
$Y^N=\{y_1,\dots,y_N\};\;y_i=x_i-\theta_N$. Then we divide the sample $Y^N$ into two sub-samples $I_1(b),\,I_2(b)$ as follows:
$$
y_i\in \left\{
\begin{array}{ll}
& I_1(b)=\{\tilde{y}_1,\dots,\tilde{y}_{N_1(b)}\},\quad -\phi(b)\le y_i\le b \\
& I_2(b)=\{\hat{y}_1,\dots,\hat{y}_{N_2(b)}\},\quad y_i>b \text{ or } y_i<-\phi(b),
\end{array}
\right.
$$
where the function $\phi(b)$ is defined from the following condition: $0=\intl_{-\phi(b)}^b\,y\,f_0(y)dy$, $f_0(y)=f_0(x-\theta_N)$, $N=N_1(b)+N_2(b)$ and $N_1(b),N_2(b)$ are sample sizes of $I_1(b),\,I_2(b)$, respectively.

2. As before we compute the statistic
$$
\Psi_N(b)=\dpfrac 1{N^2}(N_2(b) \suml_{i=1}^{N_1(b)}\,\tilde y_i-N_1(b)\suml_{i=1}^{N_2(b)}\,\hat y_i).
$$

3. Then the value $J=\max_b |\Psi_N(b)|$ is compared with the threshold $C$. If $J\le C$ then the hypothesis $H_0$ (no abnormal observations) is
accepted; if, however, $J>C$ then the hypothesis $H_0$ is rejected and the estimate of the parameter $\epsilon$ is constructed.

4. For this purpose, define the value $b_N^*$:
$$
b_N^* \in \arg\max_{b>0}|\Psi_N(b)|.
$$

Then
$$
\epsilon_N^*=N_2(b_N^*)/N.
$$

\vspace{0.5cm}

Consider application of this method for the study of the classic  $\epsilon$-contamination model:
$$
f_{\epsilon}(\cdot)=(1-\epsilon){\cal N}(\mu,\sigma^2)+\epsilon {\cal N}(\mu, \Lambda^2), \quad \Lambda^2 >> \sigma^2,\quad 0\le\epsilon <1/2.
$$

For this model, the method described above has the form:

1. From the sample of observations $X^N=\{x_1,\dots,x_N\}$ the mean value estimate $\hat\mu=\sum_{i=1}^N\,x_i/ N$ was computed.

2. The sequence $y_i=(x_i-\hat\mu)^2, \;i=1,\dots,N$ and its empirical mean $\theta_N=\sum_{i=1}^N\,y_i/ N$ are computed.

3. Then for each $b\in [0,B_{max}]$, where $B_{max}$ is a certain a
priori chosen maximal value of the parameter $b$, the sample
$Y^N=\{y_1,\dots,y_N\}$ is divided into two sub-samples in the
following way: for $\theta_N (1-\phi(b))\le y_i\le \theta_N(1+b)$
put $\tilde y_i=y_i$ (the size of the sub-sample $N_1=N_1(b)$),
otherwise put $\hat y_i=y_i$ (the size of the sub-sample
$N_2=N_2(b)$). Here we choose the function $\phi(b)$ from the
following condition:
$$
\intl_{\theta_N (1-\phi(b))}^{\theta_N (1+b)}\,y f_0(y)dy=0,
$$
where $f_0(\cdot)=N(0,(1-\epsilon)^2\sigma^2)$.

From here we obtain:
$$
\phi(b)=1-\dpfrac b{e^b-1}.
$$

4. For any $b\in [0,B_{max}]$, the following statistic is computed:
$$
\Psi_N(b)=\dpfrac 1{N^2}(N_2 \suml_{i=1}^{N_1}\,\tilde y_i-N_1\suml_{i=1}^{N_2}\,\hat y_i).
$$
where $N=N_1+N_2,\,N_1=N_1(b),\,N_2=N_2(b)$ are sizes of sub-samples of ordinary and abnormal observations, respectively.

5. Then, as above, the threshold $C>0$ is chosen and compared with the value $J=\max_b |\Psi_N(b)|$. If $J\le C$ then the hypothesis $H_0$ (no
abnormal observations) is accepted; if, however, $J>C$ then the hypothesis $H_0$ is rejected and the estimate of the parameter $\epsilon$ is
constructed as follows.

Define the value $b_N^*$:
$$
b_N^* \in \arg\max_{b>0}|\Psi_N(b)|.
$$

Then
$$
\epsilon_N^*=N_2(b_N^*)/N.
$$

\textbf{Remark.} For estimation of the threshold, we use the approach described in 2.1.3.

\bigskip In experiments the critical values of the statistic $\max_b | \Psi_N(b)|$ were computed. For this purpose, as above, for homogenous
samples (for $\epsilon=0$), $\alpha$-quantiles of the decision statistic $\max_b | \Psi_N(b)|$ were computed ($\alpha=0.95,\,0.99$). The results
obtained in 5000 trials of each test are presented in table 3.

\bigskip
{\bf Table 3.}

\bigskip
\hspace{-0.8cm}
\begin{tabular}{|c|c|c|c|c|c|c|c|c|c|}
\hline
$N$ & 50 & 100 &  300 & 500 & 800 & 1000 & 1200 & 1500 & 2000 \\
\hline
$0.95$ & 0.3031 & 0.2330 & 0.1570 & 0.1419 & 0.1252 & 0.1244 & 0.1146 & 0.1107 & 0.1075 \\
\hline
$0.99$ & 0.3699 & 0.2862 & 0.1947 & 0.1543 & 0.1436 & 0.1331 & 0.1269 & 0.1190 & 0.1157 \\
\hline
\end{tabular}

\bigskip
The quantile value for $\alpha=0.95$ was chosen as the critical threshold $C$ in experiments with non-homogenous samples (for $\epsilon\ne 0$).
For different sample sizes in 5000 independent trials of each test, the estimate of type 2 error $w_2$ ( i.e. the frequency of the event
$\maxl_b | \Psi_N(b)| <C$ for $\epsilon
>0$) and the estimate $\hat\epsilon$ of the parameter $\epsilon$ were computed. The results are presented in tables 4 and 5.

\vspace{0.5cm} {\bf Table 4.}

\bigskip
\begin{tabular}{|c|c|c|c|c|}
\hline
$\Lambda=3.0 $ & \multicolumn{4}{|c|} {$\epsilon=0.05$}  \\
\hline
$N$ & 300 & 500 &  800 & 1000 \\
\hline
$C$ & 0.1570 & 0.1419 & 0.1252 & 0.1244  \\
\hline
$w_2$ & 0.27 & 0.15 & 0.06 & 0.04 \\
\hline
$\hat\epsilon$ & 0.064 & 0.056 & 0.052 & 0.05 \\
\hline
\end{tabular}

\bigskip
{\bf Table 5.}

\bigskip
\begin{tabular}{|c|c|c|c|c|c|}
\hline
$\Lambda=5.0 $ & \multicolumn{5}{|c|} {$\epsilon=0.01$}  \\
\hline
$N$ & 1000 & 1200 &  1500 & 2000 & 3000 \\
\hline
$C$ & 0.1244 & 0.1146 & 0.1107 & 0.1075 & 0.1019 \\
\hline
$w_2$ & 0.25 & 0.20 & 0.15 & 0.10 & 0.04 \\
\hline $\hat\epsilon$ & 0.0135 & 0.013 & 0.012 & 0.011 & 0.010
 \\
\hline
\end{tabular}

\vspace{1.5cm} {\large \bf 2.2. Multiple switchings}

Suppose we obtain the data $X^N=\{x_1,\dots,x_N\}$, where the d.f.
of an observation $x_i$ can be written as follows:
$$
f(x_i)=(1-\epsilon_1-\dots-\epsilon_k)\,f_0(x_i-h_1)+\epsilon_1\,f_0(x_i-h_2)+\dots+\epsilon_k\,f_0(x_i-h_{k+1}),
$$
where $\epsilon_1\ge\epsilon_2\ge \dots \ge \epsilon_k \ge 0$,
$0\le \epsilon_1+\dots+\epsilon_k <1$, $|h_1|< |h_2|< \dots
<|h_{k+1}|$.

Suppose that the d.f. $f_0(x)$ is symmetric w.r.t. $x=0$ and
$\minl_{1\le i\le k}\,(|h_{i+1}|-|h_i|)\ge B>0$.

Our goal is to test the hypothesis $\epsilon_s=0,s=1,\dots,k$ (no switches) and in case this hypothesis is rejected to estimate the number of
switches $k\ge 1$ and the parameters of the model $\epsilon_i,\,i=1,\dots,k$, and $h_j,\,j=1,\dots,k+1$. In this section we denote by
$\mathbf{E}_i,\,i=0,1,\dots,k$, the mathematical expectation of random variables corresponding to the d.f. with shift $h_i (h_0\df 0)$.

This model has the following sense. In the case of a binary switching we have ordinary and abnormal observations. In the case of multiple
switchings abnormal observations are from different classes. The idea to use the sample mean as a reference point of the above described method
is no more valid, because in case of many classes it can be greatly biased towards the maximal $|h_i|$. Instead, we use the reference points
from the histogram of the obtained sample. Concretely, we do as follows.

1. Construct the histogram $hist_N(t)$ of data by the whole sample
$X^N$ obtained. Find $\arg\maxl_t\,hist_N(t)$. An arbitrary point
from this set is assumed to be the reference point $\theta_N$ used
in the following algorithm for a binary switching model.

1.1.Fix the parameter $b>0$ and classify observations as follows:
if an observation falls into the interval
$(\theta_N-b,\theta_N+b)$, then we place it into the sub-sample of
ordinary observations, otherwise - to the sub-sample of abnormal
observations.

1.2. Then for each $b>0$ we obtain the following decomposition of
the sample $X^N$ into two sub-samples
$$\begin{array}{ll}
& X_1(b)=\{\tilde x_1,\tilde x_2,\dots,\tilde x_{N_1}\},\quad |\tilde x_i-\theta_N|< b, \\
& X_2(b)=\{\hat x_1,\hat x_2,\dots,\hat x_{N_2}\},\quad |\hat
x_i-\theta_N|\ge b
\end{array}
$$
Denote by $N_1=N_1(b),\,N_2=N_2(b),\,N=N_1+N_2$ the sizes of the
sub-samples $X_1$ and $X_2$, respectively.

The parameter  $b$  is chosen so that the sub-samples $X_1$ and
$X_2$ are separated in the best way. For this purpose, consider
the following statistic:
$$
\Psi_N(b)=\dpfrac 1{N^2}(N_2 \suml_{i=1}^{N_1}\,\tilde
x_i-N_1\suml_{i=1}^{N_2}\,\hat x_i).
$$

1.3. Define the boundary $C>0$ and compare it with the value
$J=\max |\Psi_N(b)|$ on the set $0<b\le B$. If $J\le C$ then we
accept the hypothesis $H_0$ about the absence of abnormal
observations; if, however, $J>C$ then the hypothesis $H_0$ is
rejected and the estimates of the parameters
$\epsilon\df(\epsilon_1+\dots+\epsilon_k)$ and $h_1$ are constructed. Remark
that our primary goal is to separate ordinary and abnormal
observations in the sample. Evidently, classification errors must
be small and therefore we have to require some kind of convergence
of the estimate $\hat\epsilon_N$ to its true value
$\epsilon$.

1.4. Define the number $b_N^*$:
$$
b_N^* \in \arg\max_{0<b\le B}|\Psi_N(b)|.
$$

Then
$$
\hat{\epsilon}_N=N_2(b_N^*)/N.
$$

2. As a result, we obtain two classes of observations at the first
step (ordinary and abnormal data) and the estimate
$\hat\epsilon_N$ of the sum $\epsilon=(\epsilon_1+\dots+\epsilon_k)$, as well
as the estimate of the average $\mathbf{E}_0\,x_i$.

3. Then we remove all found 'ordinary' observations from the
sample and repeat steps 1 and 2. As a result, we obtain the
estimate $\hat\epsilon_1$ of the parameter $\epsilon_1$, as well
as the estimate of the average $\mathbf{E}_1\,x_i$.

4. So we proceed further until a sub-sample without switches is
obtained (i.e. the decision threshold $C$ is not exceeded). As a
result, we obtain the estimate $\hat k_N$ of the number of classes
$k$, as well as the estimates of the parameters $\epsilon_1>\dots
>\epsilon_k >0$ and averages $\mathbf{E}_0\,x_i,\,\mathbf{E}_1\,x_i,\dots,\mathbf{E}_k\,x_i$.

We see that this method is based upon reduction to the case of a binary switching model. In this case we characterize the quality of a method by
the performance criteria of the right estimation of the number of classes (i.e. $\hat k_N=k$) and the accuracy of estimation (e.g.,
$\maxl_i\,|\hat\epsilon_i-\epsilon_i|$ in the case $\hat k_N=k$). So we must use  the following performance criterion:
$$
\mathbf{P}_{\epsilon}\{(\hat k_N\ne k)\cup
\left((\maxl_i\,|\hat\epsilon_i-\epsilon_i|>\delta) \cap (\hat k_N=k)\right)\}.
$$

However, we see that the crucial thing is to correctly estimate the number of classes $k$. The estimates of the parameters $h_1,\dots,h_{k+1}$
are assumed to be the reference points at each step of the above described recurrent procedure. Then consistent estimates of $\epsilon_i$ can be
obtained by some of standard methods (e.g., the method of moments). Therefore we use the following performance criterion:
$$
\mathbf{P}_{\epsilon}\{\hat k_N\ne k\}.
$$

Remark that the 1st type error for multiple switchings can be estimated like in the binary case (we do not formulate this result). As to the 2nd
type error (i.e. the probability that we stop at the 1st step of the method because the decision threshold is not exceeded) just observe that a
binary switch is a particular case of the general multiple switching situation (when all $\epsilon_i$ beginning from $i=2$ are equal to zero).

Therefore
$$\begin{array}{ll}
\mathbf{P}_{\epsilon}\{2nd \text{ type error, multiple switches}\}\le & \mathbf{P}_{\epsilon}\{2nd \text{ type error, binary case}\}\\
&\le L_1 \exp(-L(\delta)N),
\end{array}
$$
for $0\le\delta\le \maxl_{0\le b\le B_{max}}\,|\Psi(b)|-C$, where
$L(\delta)=min(\dpfrac{\delta^2}{16\phi_0^2g},\dpfrac
{H\delta}{8\phi_0}),\;L_1=4\phi_0$.

Now consider the event $\{\hat k_N \ne k\}=\{k_N< k\}\cup \{k_N>
k\}$.

The event $\{k_N < k\}$ means that at a certain recurrent step of
the above described procedure a sub-sample of remaining
observations (after eliminations of all previous sub-samples) is
considered to be "pure" (i.e. without switches) but in reality it
contains some more switches. The probability of this event is less
than the 2nd type error at this step of the procedure. Therefore,
$$
\mathbf{P}_{\epsilon}\{k_N< k\}\le L_1 \exp(-L(\delta)N),
$$
for $0\le\delta = \maxl_{0\le b\le B_{max}}\,|\Psi(b)|-C$, where
$L(\delta)=min(\dpfrac{\delta^2}{16\phi_0^2g},\dpfrac
{H\delta}{8\phi_0}),\;L_1=4\phi_0$.

The event $\{k_N >k\}$ means that finally some more switches are
detected in the obtained sample than in reality. The probability
of this event is less than the 1st type error at the final step of
the above recurrent procedure:
$$
\mathbf{P}_{\epsilon}\{k_N> k\}\le L_1 \exp(-L(C)N),
$$
where $L(C)=min(\dpfrac{C^2}{16\phi_0^2g},\dpfrac
{HC}{8\phi_0}),\;L_1=4\phi_0$.

Therefore the following theorem holds.

\bigskip
{\bf Theorem 4.}

Suppose $0<C<\maxl_{0\le b\le B_{max}}\,|\Psi(b)|$. Then the 2nd
type error probability is estimated from above as follows:
$$
\mathbf{P}_{\epsilon}\{\text { 2nd type error }\}\le L_1 \exp(-L(\delta)N),
$$
where $0\le\delta = \maxl_{0\le b\le B_{max}}\,|\Psi(b)|-C$,
$L(\delta)=min(\dpfrac{\delta^2}{16\phi_0^2g},\dpfrac
{H\delta}{8\phi_0}),\;L_1=4\phi_0$.

Moreover, the estimate of the number of switchings $\hat k_N$
converges a.s. to the true value of $k$ as $N\to\infty$ and
$$
\mathbf{P}_{\epsilon}\{k_N\ne k\}\le L_1 (\exp(-L(\delta) N)+\exp(-L(C)N)),
$$
where $0\le\delta = \maxl_{0\le b\le B_{max}}\,|\Psi(b)|-C$ and
$L(C)=min(\dpfrac{C^2}{16\phi_0^2g},\dpfrac
{HC}{8\phi_0}),\,L(\delta)=min(\dpfrac{\delta^2}{16\phi_0^2g},\dpfrac
{H\delta}{8\phi_0}),\;L_1=4\phi_0$.

\bigskip{\bf Example}

Let us consider the following example. Suppose we have the model
with three classes of observations:
$$
f(x_i)=(1-\epsilon_1-\epsilon_2)\,f_0(x_i-h_1)+\epsilon_1\,f_0(x_i-h_2)+\epsilon_2\,f_0(x_i-h_3),
\qquad i=1,\dots,N,
$$
where $f_0(\cdot)={\cal N}(0,1)$; $x_i$ are i.r.v.'s.

The problem is to estimate the unknown number of classes $k=3$,
parameters $h_1,h_2,h_3$, and $\epsilon_1,\epsilon_2$ by the
sample $X^N=\{x_1,\dots,x_N\}$.

Concretely, in this model the following parameters were chosen:
$$\begin{array}{ll}
& \epsilon_1=0.3;\;\epsilon_2=0.15 \\
& h_1=1,\;h_2=3,\;h_3=7.
\end{array}
$$

For estimation of the decision threshold, the above empirical formula (2) can be used:
$$
log(C)=-0.9490-0.4729*log(N)+1.0627*log(\sigma)-0.6502*log(1-\rho)-0.2545*log(1-\alpha).
$$

Again remark that the elasticity coefficient for the factor $N$ is close to its theoretical value $-0.5$.

In experiments we estimated the number of classes $\hat k_N$ and
the corresponding error $\hat er_N=\mathbf{P}_{\epsilon}\{\hat k_N\ne k\}$.

The following results were obtained (each cell of this table is
the average in 1000 replications):

\bigskip
{\bf Table 6.}

\bigskip
\begin{tabular}{|c|c|c|c|c|c|c|c|}
  \hline
  $N$ & 100 & 200 & 300 & 500 & 700 & 1000 & 1500 \\
  \hline
  $\hat er_N$ & 0.116 & 0.090 & 0.070 & 0.048 & 0.036 & 0.016 & 0.010 \\
  \hline
\end{tabular}

\pagebreak {\large \bf 3. Multivariate models}

\bigskip
{\bf 3.1. Multivariate classification}

\bigskip
{\bf Binary mixtures}

Now let us consider the multivariate classification problem with binary mixtures. Suppose multivariate observations are of the following type:
$$
\mathbf{X}^N=\{X^n\}_{n=1}^N,\,\,X^n=(x_n^1,\dots,x_n^k).
$$

The multivariate density function of the vector $X^n$ is
$$
f(X^n)=(1-\epsilon)f_0(X^n)+\epsilon f_1(X^n),
$$
where $f_0(\cdot),\,f_1(\cdot)$ are the d.f.'s of ordinary and
abnormal observations, respectively; the d.f. $f_0(\cdot)$ is
supposed to be symmetric w.r.t. its mean vector.

First, let us consider the case $\mathbf{E}_1(X^n)=a\ne 0$, i.e. changes in mean of abnormal observations. Remark that the baseline "change-in-mean"
problem is usually considered in many methods of 'cluster analysis' in which different distances between multivariate 'points' of
characteristics (even without references to density functions and mathematical expectations of observations) are considered.

The method can be formulated in analogy with the univariate case:

1) From the initial sample $\mathbf{X}^N$  compute the estimate of the mean
value:
$$
\theta_N=\dpfrac 1N \suml_{i=1}^N\,X^i.
$$

2) Fix the parameter $b>0$ and classify observations as follows:

if $\|X^i-\theta_N\| \le b$, then we place $X^i$ into the
sub-sample of ordinary observations $\{\tilde Y^i\}$;

if $\|X^i-\theta_N\| > b$, then we place $X^i$ into the sub-sample
of abnormal observations $\{\hat Y^i\}$.

As a result, for each $b>0$ we obtain the decomposition of the
sample $\mathbf{X}^N$ into sub-samples of ordinary and abnormal
observations. Suppose the size of ordinary sub-sample is $N_1(b)$
and the size of abnormal sub-sample is $N_2(b)$.

3) The parameter  $b$  can be chosen in order to separate the
sub-samples of ordinary and abnormal observations ($\{\tilde{Y}^i\}$ and
$\{\hat{Y}^i\}$, respectively) in the best way. For this purpose, consider
the following statistic:
$$
\Psi_N(b)=\dpfrac 1{N^2}(N_2 \suml_{i=1}^{N_1}\,\tilde
Y^i-N_1\suml_{i=1}^{N_2}\,\hat Y^i).
$$

4) Define the boundary $C>0$ and compare it with the value
$J=\maxl_b \|\Psi_N(b)\|$ on the set $b>0$. If $J\le C$ then we
accept the hypothesis $H_0$ about the absence of abnormal
observations; if, however, $J>C$ then the hypothesis $H_0$ is
rejected and the estimates of the parameters $\epsilon$ and $a$
are constructed.

Remark that our primary goal is to separate ordinary and abnormal
observations in the sample. Evidently, classification errors must
be small and therefore we have to require some kind of convergence
of the estimate $\hat\epsilon_N$ to its true value $\epsilon$.

5) Define the number $b_N^*$:
$$
b_N^* \in \arg\max_{b>0}\|\Psi_N(b)\|.
$$

Then
$$
\epsilon_N^*=N_2(b_N^*)/N,\quad a_N^*=\theta_N / \epsilon_N^*.
$$
are the nonparametric estimates for $\epsilon$ and $a$,
respectively.

Our main results in this case are analogous to the univariate situation.

\bigskip
{\bf Theorem 5. }

Suppose $\epsilon=0$ and the d.f. $f_0(\cdot)$ is symmetric w.r.t. its mean vector. Then for any $C>0$ the following upper estimate for the
probability of the 1st type error holds:
$$
\mathbf{P}_0\{\maxl_{b>0}\,\|\Psi_N(b)\| >C\} \le 4\phi_0(CN/2) \exp(-L(C)N),
$$
where $L(C)=min\left(\dpfrac {HC}{8\phi_0(CN/2)},\;\dpfrac {C^2}{16\phi_0^2(CN/2)g}\right)$, the constants $g,H$ are taken from the uniform
Cramer  condition.

For the 2nd type error we can formulate the following result.

\bigskip
{\bf Theorem 6.}

Suppose all assumptions of theorem 5 are satisfied and there exists $r^*=\supl_b\,r(b)$. Suppose also that $f^{''}(\cdot)\not=0$ and continuous.
Then for $0< C< \maxl_b\,|\Psi(b)|$ we have

1)
$$
\mathbf{P}_{\epsilon}\{\maxl_b\,\|\Psi_N(b)\|\le C\}\le 4\phi_0(CN/2+r^*)\,\exp(-L(\delta) N)
$$
where $\delta=\maxl_b\,\|\Psi(b)\|-C >0,\;L(\delta)=min(\dpfrac {\delta^2}{16\phi_0^2g},\dpfrac {H\delta}{8 \phi_0})$.

2) Suppose, moreover, that equation (*) has a unique root $b^{*}$. Then

 $b_N^{*}\to b^{*}$ $\mathbf{P}_{\epsilon}$-a.s. as $N\to\infty$;

This method deals with binary mixtures of multivariate d.f.'s. Its generalization to multiple classes of multivariate d.f.'s can be obtained in
analogy with the previous section.

\bigskip
{\bf Multiple switches}

In this case the multivariate density function of the vector $X^n$ is
$$
f(X^n)=(1-\epsilon_1-\dots-\epsilon_k)f_0(X^n-h_1)+\epsilon_1\,f_0(X^n-h_2)+\dots+\epsilon_k\,f_0(X^n-h_{k+1})
$$

where $\epsilon_1\ge\epsilon_2\ge \dots \ge \epsilon_k \ge 0$, $0\le \epsilon_1+\dots+\epsilon_k <1$, $\|h_1\|< \|h_2\|< \dots <\|h_{k+1}\|$.

Suppose that the d.f. $f_0(x)$ is symmetric w.r.t. $x=0$ and $\minl_{1\le i\le k}\,(\|h_{i+1}\|-\|h_i\|)\ge B>0$.

In order to estimate the number of classes $k$, as well as parameters $\epsilon_i$ we do as follows:

\bigskip

From the sample of initial multivariate observations
$$
X^n=(x_n^1,\dots,x_n^k), \qquad n=1,\dots,N.
$$
we build the sample of their Euclidean norms:
$$
Y_n=\|X^n\|=\sqrt {(x_n^1)^2+\dots+(x_n^k)^2}, \qquad n=1,\dots,N.
$$

1. Construct the histogram $hist_N(t)$ of data by the whole sample $Y^N=(Y_1,\dots,Y_N)$. Find $\arg\maxl_t\,hist_N(t)$. An arbitrary point from
this set is assumed to be the reference point $\theta_N$ used in the following algorithm for a binary switching model.

1.1. Fix the parameter $b>0$ and classify observations as follows:

if $\|Y_i-\theta_N\| \le b$, then we place $Y_i$ into the sub-sample of ordinary observations ($\tilde Y_N^i$);

if $\|Y_i-\theta_N\| > b$, then we place $Y_i$ into the sub-sample of abnormal observations ($\hat Y_N^i$).

1.2. Then for each $b>0$ we obtain the following decomposition of the sample $Y^N$ into two sub-samples
$$\begin{array}{ll}
& Y_1(b)=\{\tilde Y_1,\tilde Y_2,\dots,\tilde Y_{N_1}\},\quad \|\tilde Y_i-\theta_N\|< b, \\
& Y_2(b)=\{\hat Y_1,\hat Y_2,\dots,\hat Y_{N_2}\},\quad \|\hat Y_i-\theta_N\|\ge b
\end{array}
$$
Denote by $N_1=N_1(b),\,N_2=N_2(b),\,N=N_1+N_2$ the sizes of the sub-samples $Y_1$ and $Y_2$, respectively.

The parameter  $b$  is chosen so that the sub-samples $Y_1(b)$ and $Y_2(b)$ are separated in the best way. For this purpose, consider the
following statistic:
$$
\Psi_N(b)=\dpfrac 1{N^2}(N_2 \suml_{i=1}^{N_1}\,\tilde Y_i-N_1\suml_{i=1}^{N_2}\,\hat Y_i).
$$

1.3. Define the boundary $C>0$ and compare it with the value $J=\max |\Psi_N(b)|$ on the set $0<b\le B$. If $J\le C$ then we accept the
hypothesis $H_0$ about the absence of abnormal observations; if, however, $J>C$ then the hypothesis $H_0$ is rejected and the estimates of the
parameters $\epsilon=(\epsilon_1+\dots+\epsilon_k)$. Remark that our primary goal is to separate ordinary and abnormal observations in the
sample. Evidently, classification errors must be small and therefore we have to require some kind of convergence of the estimate
$\hat\epsilon_N$ to its true value $\epsilon_1+\dots+\epsilon_k$.

1.4. Define the number $b_N^*$:
$$
b_N^* \in \arg\max_{0<b\le B}\|\Psi_N(b)\|.
$$

Then
$$
\epsilon_N^*=N_2(b_N^*)/N.
$$

2. As a result, we obtain two classes of observations at the first step (ordinary and abnormal data) and the estimate $\hat\epsilon_N$ of the
sum $\epsilon_1+\dots+\epsilon_k$.

3. Then we remove all found 'ordinary' observations from the sample and repeat steps 1 and 2. As a result, we obtain the estimate
$\hat\epsilon_1$ of the parameter $\epsilon_1$.

4. So we proceed further until a sub-sample without switches is obtained (i.e. the decision threshold $C$ is not exceeded). As a result, we
obtain the estimate $\hat k_N$ of the number of classes $k$, as well as the estimates of the parameters $\epsilon_1>\dots
>\epsilon_k >0$.

Again we remark that the 1st type error for multiple switchings can be estimated like in the binary case (we do not formulate this result). As
to the 2nd type error (i.e. the probability that we stop at the 1st step of the method because the decision threshold is not exceeded) just
observe that a binary switch is a particular case of the general multiple switching situation (when all $\epsilon_i$ beginning from $i=2$ are
equal to zero.

Therefore
$$\begin{array}{ll}
\mathbf{P}_{\epsilon}\{2nd \text{ type error, multiple switches}\}\le & \mathbf{P}_{\epsilon}\{2nd \text{ type error, binary case}\}\\
&\le L_1 \exp(-L(\delta)N),
\end{array}
$$
for $0\le\delta\le \maxl_{0\le b\le B}\,\|\Psi(b)\|-C$, where $L(\delta)=min(\dpfrac{\delta^2}{16\phi_0^2g},\dpfrac
{H\delta}{8\phi_0}),\;L_1=4\phi_0$.

Now consider the event $\{\hat k_N \ne k\}=\{k_N< k\}\cup \{k_N> k\}$.

The event $\{k_N < k\}$ means that at a certain recurrent step of the above described procedure a sub-sample of remaining observations (after
eliminations of all previous sub-samples) is considered to be "pure" (i.e. without switches) but in reality it contains some more switches. The
probability of this event is less than the 2nd type error at this step of the procedure. Therefore,
$$
\mathbf{P}_{\epsilon}\{k_N< k\}\le L_1 \exp(-L(\delta)N),
$$
for $0\le\delta = \maxl_{0\le b\le B}\,\|\Psi(b)\|-C$, where $L(\delta)=min(\dpfrac{\delta^2}{16\phi_0^2g},\dpfrac
{H\delta}{8\phi_0}),\;L_1=4\phi_0$.

The event $\{k_N >k\}$ means that finally some more switches are detected in the obtained sample than in reality. The probability of this event
is less than the 1st type error at the final step of the above recurrent procedure:
$$
\mathbf{P}_{\epsilon}\{k_N> k\}\le L_1 \exp(-L(C)N),
$$
where $L(C)=min(\dpfrac{C^2}{16\phi_0^2g},\dpfrac {HC}{8\phi_0}),\;L_1=4\phi_0$.

Therefore the following theorem holds.

\bigskip
{\bf Theorem 7.}

Suppose $0<C<\maxl_{0\le b\le B}\,\|\Psi(b)\|$. Then the 2nd type error probability is estimated from above as follows:
$$
\mathbf{P}_{\epsilon}\{\text { 2nd type error }\}\le L_1 \exp(-L(\delta)N),
$$
where $0\le\delta = \maxl_{0\le b\le B}\,\|\Psi(b)\|-C$, $L(\delta)=min(\dpfrac{\delta^2}{16\phi_0^2g},\dpfrac
{H\delta}{8\phi_0}),\;L_1=4\phi_0$.

Moreover, the estimate of the number of switchings $\hat k_N$ converges a.s. to the true value of $k$ as $N\to\infty$ and
$$
P_{\epsilon}\{k_N\ne k\}\le L_1 (\exp(-L(\delta) N)+\exp(-L(C)N)),
$$
where $0\le\delta = \maxl_{0\le b\le B}\,\|\Psi(b)\|-C$ and $L(C)=min(\dpfrac{C^2}{16\phi_0^2g},\dpfrac
{HC}{8\phi_0}),\,L(\delta)=min(\dpfrac{\delta^2}{16\phi_0^2g},\dpfrac {H\delta}{8\phi_0}),\;L_1=4\phi_0$.

\bigskip
{\bf  Example}

 Suppose we have the model with three classes of multivariate Gaussian observations:
$$
f(x_i)=(1-\epsilon_1-\epsilon_2)\,f_0(x_i-h_1)+\epsilon_1\,f_0(x_i-h_2)+\epsilon_2\,f_0(x_i-h_3),
\qquad i=1,\dots,N,
$$
where $f_0(\cdot)$ has the multivariate Gaussian d.f. with the
vector of means $\mu=(\mu_1,\mu_2)^{'}$ and the covariance matrix
$Cov(x_i)=\left(
\begin{array}{cc}
  0.745 & -0.07 \\
  -0.07 & 0.51 \\
\end{array}
\right) $ .

The problem is to estimate the unknown number of classes $k=3$,
parameters $h_1,h_2,h_3$, and $\epsilon_1,\epsilon_2$ by the
sample $X^N=\{x_1,\dots,x_N\}$.

Concretely, in this model the following parameters were chosen:
$$\begin{array}{ll}
& \epsilon_1=0.3;\;\epsilon_2=0.15 \\
& h_1=(0\; 0)^{'},\;h_2=(1\; 2)^{'},\;h_3=(2\; 3)^{'}.
\end{array}
$$

We take the norm of the vectors $x_i$ and so reduce this problem
to the univariate case considered earlier in this paper.

For estimation of the decision threshold the above formula (2) can be used:
$$
log(C)=-0.9490-0.4729*log(N)+1.0627*log(\sigma)-0.6502*log(1-\rho)-0.2545*log(1-\alpha).
$$

Again we remark that the main problem is to estimate the number of
classes $\hat k_N$ (estimation of $h_i$ and $\epsilon_j$ can be
done with the help of some standard methods for the given model
structure).

In experiments we estimated the number of classes $\hat k_N$ and
the corresponding error $\hat er_N=P_{\epsilon}\{\hat k_N\ne k\}$.

The following results were obtained (each cell of this table is
the average in 1000 independent trials of the test):

\bigskip
{\bf Table 7.}

\bigskip
\begin{tabular}{|c|c|c|c|c|c|c|c|}
  \hline
  $N$ & 100 & 200 & 300 & 500 & 700 & 1000 & 1500 \\
  \hline
  $\hat er_N$ & 0.991 & 0.910 & 0.707 & 0.189 & 0.049 & 0.020 & 0.004 \\
  \hline
\end{tabular}

\vspace{2cm} {\large \bf 3.2. Switching regressions}

The following model of observations was considered:
$$
y_i=X\beta_i+u_i=X(\dzeta_i\beta_0+(1-\dzeta_i)\beta_1)+u_i,
$$
where

$y$  is a $N\times 1$ vector of dependent observations;

$X$  is a $N\times k$ matrix of predictors;

$u$  is a $N\times 1$ vector of centered random noises;

$\beta_i$  is a $k\times 1$ vector of model coefficients, $\dzeta_i$ is a Bernoulli distributed r.v. (independent from $u_i$) with two states:
$1$ with the probability $(1-\epsilon)$ and $0$ with the probability $\epsilon$ for a certain unknown parameter $0<\epsilon<1$. Here $\beta_0\ne
\beta_1$.

In terms, we suppose that regression coefficients of this model
can change (switch) form the level $\beta_0$ to $\beta_1$ and the
mechanism of this change is purely random. We need to test the
hypothesis about the absence of switchings for each coefficient
($\epsilon=0$) and in the case of rejection of this hypothesis to
construct the estimate of the parameter $\epsilon>0$.

For solving this problem, consider the OLS estimate of the vector $\beta_i$ (here and below $'$ is the symbol of transposition):
$$
\hat\beta_i=(X'X)^{-1}X'y_i=\dzeta_i\beta_0+(1-\dzeta_i)\beta_1+(X'X)^{-1}X'u_i.
$$

Since the sequence of noises $u$ is centered, the problem is reduced to the above considered problem of detection of switches in the mean of an
observed random vector. The matrix of predictors $X$ influences only the random component.

Formally, we need to introduce the following vector
$I=(1,1,\dots,1)$ ($N$ units) and consider
$$
\tilde\beta_i=[\dzeta_i\beta_0+(1-\dzeta_i)\beta_1]\,I+(X^{'}X)^{-1}X^{'}u_i\,I.
$$

Then the $(k\times n)$ matrix $\tilde \beta_i$ consists of $N$ columns of $k\times 1$ vectors with  means $\beta_0$ and $\beta_1$ changing in a
purely random manner. Each component $j=1,\dots,k$ of these vectors $\tilde \beta_i^j,\;i=1,\dots,N$ is therefore a univariate random sequence
$$
\tilde\beta_i^j=[\dzeta_i\beta_0^j+(1-\dzeta_i)\beta_1^j]_i+\xi_i^j, \qquad i=1,\dots,N,
$$
where
$$
\xi_i^j=((X^{'}X)^{-1}X^{'}u\,I)_i^j.
$$

So the problem of detection of changes in regression coefficients is reduced to the above considered problem of detection switches in the mean
value of a univariate random sequence. Remark that the uniform Cramer and the $\psi$-mixing conditions are still satisfied for the process
$\xi_i^j,\,i=1,\dots,N$. As $\mathbf{E}u_i\equiv 0$ we get that there exist constants $g_1>0,\,H_1>0$ such that
$$
Ee^{t\,\xi_i^j}\le e^{\dpfrac 12 g_1t^2},\quad |t|\le H_1,
$$
for all $i=1,\dots,N,\,j=1,\dots,k$. Moreover, we choose the
number $m_0(\cdot)$ from the $\psi$-mixing condition for
$\xi_i^j,\,i=1,\dots,N$: for any chosen number $\gamma(x)>0$:
$\psi(l)\le \gamma(x)$ for $l\ge m_0(x)$.

For testing the hypothesis of no switches we again consider the
decision statistic $\Psi_N(b)$ and compare the maximum of its
module with the decision threshold $C>0$. Then the following
theorem holds:

\bigskip
{\bf Theorem 8.}

Suppose $\epsilon=0$, the d.f. of $u_i$ is symmetric w.r.t. zero
anf the $\psi$-mixing and the uniform Cramer conditions for
$\xi_i^j,\,i=1,\dots,N$ are satisfied. Then for any threshold
$C>0$ the following upper estimate for the 1st type error
probability holds:
$$
\mathbf{P}_0\{\maxl_{b>0}\,|\Psi_N(b)| >C\} \le 4m_0(CN/2)
\exp(-L(C)N),
$$
where $L(C)=min\left(\dpfrac {H_1C}{8m_0(CN/2)},\;\dpfrac
{C^2}{16m_0^2(CN/2)g_1}\right)$, the constants $g_1,H_1$ are taken
from the uniform Cramer  condition.

\bigskip
In order to consider the 2nd type error we just remark that the considered switching regression model is equivalent to the following
specification of a model with the binary switches in mean:
$$
f_{\tilde\beta_i^j}(x)=(1-\epsilon)f_{\xi_i^j}(x-\beta_0^j)+\epsilon\,f_{\xi_i^j}(x-\beta_1^j).
$$

Denote $h^j=\beta_1^j-\beta_0^j\ne 0$ and consider the value
$$
r_{\tilde\beta_i^j}(b)=\intl_{\beta_0^j+\epsilon
h^j-b}^{\beta_0^j+\epsilon h^j+b}\,f_{\tilde\beta_i^j}(x)x\,dx.
$$

Then the following theorem holds.

\bigskip
{\bf Theorem 9.}

 Suppose all assumptions of theorem 8 are satisfied and there
exists $r^*=\supl_b\,r_{\tilde\beta_i^j}(b)$. Suppose also that
$f_{\xi_i^j}^{''}(\cdot)\not=0$ and continuous. Then for $0< C<
\maxl_b\,|\Psi_{\tilde\beta_i^j}(b)|$ we have
$$
\mathbf{P}_{\epsilon}\{\maxl_b\,|\Psi_N(b)|\le C\}\le
4\phi_0(CN/2+r^*)\,\exp(-L(\delta) N)
$$
where $\delta=\maxl_b\,|\Psi_{\tilde\beta_i^j}(b)|-C
>0,\;L(\delta)=min(\dpfrac {\delta^2}{16m_0^2g_1},\dpfrac
{H_1\delta}{8 m_0})$.

\bigskip
{\bf Example}

In the following example the regression model with one
deterministic predictor was considered:
$$
y_i=c_1+c_2*i+u_i, \quad u_i\sim N(0;1),\quad i=1,\dots,n.
$$

$$
\xi\sim U[0;1]
$$
$$
\beta=[c_1;\;c_2]=\left\{
\begin{array}{ll}
& \beta_1,\quad \epsilon_1<\xi\le 1 \\
& \beta_2,\quad 0\le\xi\le\epsilon_1
\end{array}
\right.
$$

 {\bf Table 8.}

\bigskip
\begin{tabular}{|c|c|c|c|c|}
\hline
$\epsilon=0.05 $ & \multicolumn{4}{|c|} {$\beta_1=[1; 1],\;\beta_2=[1; 2]$}  \\
\hline
$N$ & 300 & 500 &  800 & 1000 \\
\hline
$C$ & 0.07 & 0.05 & 0.04 & 0.03  \\
\hline
$w_2$ & 0.87 & 0.59 & 0.14 & 0.004 \\
\hline
$\hat\epsilon$ & 0.08 & 0.059 & 0.052 & 0.05 \\
\hline
\end{tabular}

\bigskip
{\bf Table 9.}

\bigskip
\begin{tabular}{|c|c|c|c|c|}
\hline
$\epsilon=0.1 $ & \multicolumn{4}{|c|} {$\beta_1=[1; 1],\;\beta_2=[1; 1.5]$}  \\
\hline
$N$ & 300 & 500 &  800 & 1000 \\
\hline
$C$ & 0.07 & 0.05 & 0.04 & 0.03  \\
\hline
$w_2$ & 0.83 & 0.65 & 0.13 & 0.0 \\
\hline
$\hat\epsilon$ & 0.15 & 0.12 & 0.102 & 0.10 \\
\hline
\end{tabular}

\vspace{1.5cm}
 {\large \bf  Conclusion}

In this paper problems of the  retrospective detection/estimation of 'abnormal' observations were considered. The detection/estimation method was proposed.
It was proved that type 1 and type 2 errors of the proposed method converge to zero exponentially as the sample size $N$ tends to infinity.
The asymptotic optimality of the proposed method follows from theorem 3. In this theorem the theoretical lower bound for the error of
estimation of the model's parameters was established. This bound is attained for the proposed method (by the order of convergence to zero of
the estimation error).

\newpage

 {\bf Proofs of theorems.}

\bigskip
{\bf Proof of theorem 1.}

First, let us prove the following inequality:
$$
\maxl_{b>0}\,\mathbf{P}_0\{|\Psi_N(b)| >C \} \le L_1 \exp (-L_2(C)N),
$$
where $L_1,L_2(C)$ are some positive constant and function not depending on $N$.

For the statistic $\Psi_N(b)$ we can write:
$$
\Psi_N(b)=(N\,\suml_{i=1}^{N_1}\,\tilde x_i-N_1\suml_{i=1}^N\,x_i)/ N^2.
$$

Then
$$
\mathbf{P}_0\{|\Psi_N(b)| >C \} \le \mathbf{P}_0\{| \suml_{i=1}^{N_1}\,\tilde x_i|
>\dpfrac C2 N\}+\mathbf{P}_0\{N_1 |\suml_{i=1}^N\,x_i| >\dpfrac C2 N^2 \}.
$$

Further,
$$
\mathbf{P}_0\{| \suml_{i=1}^{N_1}\,\tilde x_i|
>\dpfrac C2 N\}=\mathbf{P}_0\{\suml_{i=1}^{N_1}\,\tilde x_i
>\dpfrac C2 N\}+\mathbf{P}_0 \{\suml_{i=1}^{N_1}\,\tilde x_i
< -\dpfrac C2 N\}.
$$

For any x>0, let us choose the number $\gamma(x)$ from the following condition:
$$
\ln(1+\gamma(x))=\left \{
\begin{array}{ll}
& \dpfrac {x^2}{4g}, \qquad x\le gH \\
& \dpfrac {xH}4, \qquad x>gH,
\end{array}
\right.
$$
where $g,H$ are taken from the uniform Cramer condition.

For the chosen $\gamma(x)$, let us find such integer $\phi_0(x)\ge 1$ from the $\psi$-mixing condition that $\psi(l)\le\gamma (x)$ for $l\ge
\phi_0(x)$. Take $x=CN/2$ and denote $\phi_0(CN/2)=\phi_0(\cdot),\,\,\gamma(CN/2)=\gamma(\cdot)$.

For some fixed $n$ denote $S_n=\suml_{i=1}^n\,\tilde x_i$ and estimate the probability $\mathbf{P}_0\{S_{n}> CN/2\}$.

Consider the following decomposition of $S_n$ into groups of weakly dependent terms:
$$\begin{array}{ll}
& S_n=S_n^{(1)}+S_n^{(2)}+\dots+S_n^{(\phi_0(\cdot))} \\
& S_n^{(i)}=\tilde x_i+\tilde x_{i+\phi_0}+\dots+\tilde x_{i+\phi_0(\cdot)[\frac {n-i}{\phi_0(\cdot)}]}, \qquad i=1,2,\dots,\phi_0(\cdot).
\end{array}
$$

The number of terms within each group is no less than $[n/ \phi_0(\cdot)]$ and no more than $[n/ \phi_0(\cdot)]+1$ and the $\psi$-mixing coefficient between
terms within each group is no more than $\gamma(\cdot)$.

Then
$$\begin{array}{ll}
& \mathbf{P}_0\{S_{n}> \dpfrac C2 N\}\le \suml_{i=1}^{\phi_0(\cdot)}\,\mathbf{P}_0\{S_{n}^{(i)}> \dpfrac{CN}{2\phi_0(\cdot)}\} \\
& \le \phi_0(\cdot) \maxl_{1\le i\le \phi_0(\cdot)}\,\mathbf{P}_0\{|S_{n}^{(i)}| \ge \dpfrac {CN}{2\phi_0(\cdot)} \}.
\end{array}
$$

Consider $Z_k^{(i)}\df\suml_{j=0}^k\,\tilde x (i+\phi_0(\cdot) j)$ and obtain the exponential upper estimate for $\mathbf{P}_0\{Z_k^{(i)} >x\},\;\forall x>0$.

In virtue of Chebyshev's inequality, we have
$$
\mathbf{P}_0\{Z_k^{(i)} >x\}\le e^{-tx}\cdot \mathbf{E}_0 e^{tZ_k^{(i)}}, \qquad \forall t>0.
$$

From $\psi$-mixing condition (see Ibragimov, Linnik, 1971) and choosing $\gamma(\cdot)$ we have
$$
\mathbf{E}_0^{tZ_k^{(i)}}\le \left(1+\gamma(\cdot)\right)^k\,\mathbf{E}_0\exp(t\tilde x(i))E\exp(t\tilde x(i+\phi_0))\dots \mathbf{E}_0\exp(t\tilde x(i+\phi_0 k)).
$$

Therefore,  for $0\le t\le H$
$$
\mathbf{E}_0\exp(tS_{n})\le \left(1+\gamma(\cdot)\right)^N\,\exp(\dpfrac 12 t^2gN).
$$

Hence,
$$
\mathbf{P}_0\{S_{n}(x)> \dpfrac C2 N\}\le \phi_0(\cdot) \left(1+\gamma(\cdot)\right)^N\,\exp\left(\dpfrac N2(t^2g-Ct/ \phi_0(\cdot))\right).
$$

Taking the maximum of the right hand w.r.t. $0\le t\le H$ and taking into account the choice of $\gamma(\cdot)$ we have
$$
\mathbf{P}_0\{\suml_{i=1}^{n}\,\tilde x_i >\dpfrac C2 N\}\le \phi_0\left\{
\begin{array}{ll}
& \exp(-\dpfrac {C^2 N}{16\phi_0^2(\cdot) g}),\qquad 0<t<gH, \\
& \exp(-\dpfrac {CHN}{8\phi_0(\cdot)}), \qquad t>gH
\end{array}
\right.
$$
As this estimate does not depend of $n$, we get
$$
\maxl_{b>0}\,\mathbf{P}_0\{|\Psi_N(b)|>C\} \le 4\phi_0(\cdot)\exp(-L(C)N),
$$
where
$$
L(C)=\min \left(\dpfrac {HC}{8\phi_0(\cdot)},\dpfrac {C^2}{16\phi^2_0(\cdot)g}\right).
$$

Note that we obtained the uniform (w.r.t. the parameter $b$) exponential upper estimate for the first type error. Therefore, the same upper
estimate is valid for the probability:
$$
\mathbf{P}_0\{\maxl_{b>0}|\Psi_N(b)| >C\}\le \mathbf{P}_0\{\maxl_{b>0}\,|\suml_{i=1}^{N_1}\,\tilde x_i| > \dpfrac C2 N\}+\mathbf{P}_0\{|\suml_{i=1}^N\,x_i| > \dpfrac C2 N\}.
$$

In fact, consider the r.v. $U_N(\omega)=\maxl_{b>0}\,|\suml_{i=1}^{N_1}\,\tilde x_i|$ and define
$$
\tau_N(\omega)=\min \{1\le n\le N: |\suml_{i=1}^n\,\tilde x_i|=U_N\}.
$$

Then
$$\begin{array}{ll}
& \mathbf{P}_0\{U_N> CN/2\}=\mathbf{P}_0\{\suml_{k=1}^N\,|\suml_{i=1}^k\,\tilde x_i|\,I\{\tau_N=k\} > CN/2\}\\
& \le \mathbf{P}_0\{|\suml_{i=1}^{k_{max}(\omega)}\,|\tilde x_i| > CN/2\},
\end{array}
$$
where for any $\omega\in\Omega$: $|\suml_{i=1}^{k_{max}(\omega)}\,\tilde x_i| =\maxl_{1\le i\le N}\,|\suml_{i=1}^k\,\tilde x_i|$.

As above, we obtain the uniform upper estimate for the probability $\mathbf{P}_0\{|\suml_{i=1}^{k_{max}(\omega)}\, \tilde x_i| > CN/2\}$. Therefore,
$$
\mathbf{P}_0\{|\maxl_{b>0}\,\Psi_N(b)|>C\} \le 4\phi_0(\cdot)\exp(-L(C)N),
$$
where
$$
L(C)=\min \left(\dpfrac {HC}{8\phi_0(\cdot)},\dpfrac {C^2}{16\phi_0^2(\cdot)g}\right).
$$

Theorem 1 is proved.

\bigskip
{\bf Proof of theorem 2.}

Consider the main statistic:
$$
\Psi_N(b)=\left(N\suml_{i=1}^{N_1(b)}\,\tilde x_i-N_1(b)\,\suml_{i=1}^N\,x_i\right)/ N^2.
$$

We have
$$
\begin{array}{ll}
& \dpfrac 1N\,\mathbf{E}_{\epsilon}\suml_{i=1}^{N_1}\,\tilde x_i=\dpfrac 1N\,\suml_{n=1}^N\,\mathbf{E}_{\epsilon}(\suml_{i=1}^n\,\tilde x_i |N_1=n)\,\mathbf{P}_{\epsilon}\{N_1=n\}=\\[3mm]
&=\dpfrac
1N\,\suml_{n=1}^N\,\suml_{i=1}^n\,\mathbf{E}_{\epsilon}(\tilde x_i
| \theta_N-b <\tilde x_i <
\theta_N+b)\mathbf{P}_{\epsilon}\{N_1=n\}=\dpfrac 1N\,\left(\mathbf{E}_{\epsilon}N_1\right)\,\intl_{\epsilon h-b}^{\epsilon h+b}\,f(x)xdx / \intl_{\epsilon h-b}^{\epsilon h+b}\,f(x)dx=\\[3mm]
&\to\intl_{\epsilon h-b}^{\epsilon h+b}\,f(x)xdx, \qquad \text{
as } N\to\infty
\end{array}
$$
Here we used the relation
$$
\dpfrac 1N \mathbf{E}_{\epsilon}N_1=\dpfrac
1N\,\mathbf{E}_{\epsilon}\suml_{k=1}^{N}k\mathbb{I}(|x_k-\theta_N|\le
b)\to \intl_{\epsilon h-b}^{\epsilon h+b}f(x)dx \qquad \text{  as
} N\to\infty
$$
Therefore, using the latter relations, taking into account the law of large numbers and the relation
$$
\mathbf{E}_{\epsilon}x_i=\epsilon h
$$
we have
$$
\mathbf{E}_{\epsilon}\Psi_N(b)\to \Psi(b) \qquad \text{  as }
N\to\infty,
$$
where $\Psi(b)=r(b)-\epsilon h\,d(b)$.

For any $C>0$ we can write:
$$
\mathbf{P}_{\epsilon}\{ |\Psi_N(b)-\Psi(b)| >C\} \le  \mathbf{P}_{\epsilon}\{ |\suml_{i=1}^{N_1(b)}\,\tilde x_i-Nr(b)| > \dpfrac C2 N\}+\mathbf{P}_{\epsilon}\{| \dpfrac
{N_1(b)}N\,\suml_{i=1}^N\,x_i-N\epsilon h d(b)| > \dpfrac C2 N\}.  \eqno (3)
$$

Consider the first term in the right hand:
$$
\mathbf{P}_{\epsilon}\{ |\suml_{i=1}^{N_1(b)}\,\tilde x_i-Nr(b)| > \dpfrac C2 N\}=\mathbf{P}_{\epsilon}\{ \suml_{i=1}^{N_1(b)}\,\tilde x_i > \dpfrac C2 N+Nr(b)\}+\mathbf{P}_{\epsilon} \{
\suml_{i=1}^{N_1(b)}\,\tilde x_i < -\dpfrac C2 N+Nr(b)\}. \eqno (4)
$$

Analogously theorem 1, we put $x=CN/2+r(b)N$, find $\phi_o(x)\df\phi_o(\cdot)$ corresponding to this $x$, decompose the sum $\suml_{i=1}^{N_1}\,\tilde x_i$ into
$\phi_0(\cdot)$ groups of weakly dependent components and for each of these groups use Chebyshev's inequality.

Using considerations analogous to those in theorem 1, finally, for large enough $N$ we obtain:
$$
\mathbf{P}_{\epsilon}\{\suml_{i=1}^{N_1(b)}\,\tilde x_i >\dpfrac C2 N+Nr(b)\}\le \phi_0(\cdot)\left\{
\begin{array}{ll}
& \exp(-\dpfrac {C^2 N}{16\phi_0^2(\cdot) g}),\qquad 0<t<gH, \\
& \exp(-\dpfrac {CHN}{8\phi_0(\cdot)}), \qquad t>gH
\end{array}
\right.
$$

The second term in the right hand of (4) is estimated from above
in the same way.

As to the second term in the right hand of (3), since $N_1(b)\le N$ for any $\omega$, we obtain an analogous exponential upper estimate for it.

Again remark that we obtained the uniform (w.r.t. $b>0$) exponential upper estimate for the error probability. Therefore as in theorem 1 we can
prove the following exponential estimate:
$$
\mathbf{P}_{\epsilon}\{\maxl_b\,|\Psi_N(b)-\Psi(b)| >C\}\le 4\phi_0(\cdot)\left\{
\begin{array}{ll}
& \exp(-\dpfrac {C^2 N}{16\phi_0^2(\cdot) g}),\qquad 0<C<gH, \\
& \exp(-\dpfrac {CHN}{8\phi_0(\cdot)}), \qquad C
>gH
\end{array}
\right.
$$

For type 2 error we can write:
$$\begin{array}{ll}
 \mathbf{P}_{\epsilon}\{\maxl_b\,|\Psi_N(b)| < C\}\le
&\mathbf{P}_{\epsilon}\{\maxl_b\,|\Psi_N(b)-\Psi(b)|
> \maxl_b|\Psi(b)|-C\}\\
&\le 4\phi_0(\cdot)\,\left \{
\begin{array}{ll}
& \exp(-\dpfrac N{16}\,\dpfrac {\delta^2}{\phi_0^2(\cdot)g)}), \quad 0<\delta\le
g H,\\
& \exp(-\dpfrac {N H \delta}{8\phi_0(\cdot)}), \quad \delta>g H,
\end{array}
\right.
\end{array}
$$
where $\delta=\maxl_b|\Psi(b)|-C$.

This completes the proof of 1).

As to the proof of 2), remark that the function $\Psi(b)=\mathbf{E}_{\epsilon}\Psi_N(b)$ satisfies the reversed Lipschitz condition in a neighborhood of $b^*$.

In
fact, we have $\Psi(b^*)=0,\,\Psi^{'}(b^*)=0$ and $\Psi^{''}(b^*)=(f(\epsilon h+b^*)-f(\epsilon h-b^*))+b^*(f^{'}(\epsilon h+b^*)-f^{'}(\epsilon
h-b^*))=2(b^*)^2\,f^{''}(u)\ne 0$, where $0\le u=u(b^*)\le b^*$. Therefore in a small neighborhood of $b^*$ we obtain:
$$
|\Psi(b)-\Psi(b^*)|=(b^*)^2\,|f^{''}(u(b^*))|(b-b^*)^2\ge C(b-b^*)^2,
$$
for a certain $C=C(b^*)>0$.

Now for any $0< \kappa <1$ consider the event $|b_N-b^*|>\kappa$. Then
$$
\mathbf{P}_{\epsilon}\{|b_N-b^*|>\kappa\}
 \le \mathbf{P}_{\epsilon}\{\maxl_b\,|\Psi_N(b_N)-\Psi(b^*)| >\dpfrac 12\,C\kappa^2\}\\
\le 4\,\phi_0(\cdot)\,\exp(-L(C)N),
$$
where $L(C)=min(\dpfrac {C^2 \kappa^4}{8\phi_0^2(\cdot)g},\dpfrac {HC\kappa^2}{16\phi_0(\cdot) })$.

 From this inequality it follows that $b_N\to b^*$ $\mathbf{P}_{\epsilon}-$a.s. as $N\to\infty$.

Then
$$
\epsilon_N=N_2(b_N)/N, \qquad h_N=\theta_N/ \epsilon_N
$$
are the nonparametric estimates for $\epsilon$ and $h$, respectively.

In general these estimates are asymptotically biased and non-consistent. For construction of consistent estimates of $\epsilon$ and $h$, we need
information about the d.f. $f_0(\cdot)$. These consistent estimates can be obtained from the following system of equations:
$$\begin{array}{ll}
& \hat\epsilon_N \hat h_N=\theta_N \\
& \dpfrac{1-\hat\epsilon_N}{\hat\epsilon_N}=\dpfrac {f_0(\theta_N-b_N-\hat h_N)-f_0(\theta_N+b_N-\hat
h_N)}{f_0(\theta_N+b_N)-f_0(\theta_N-b_N)}.
\end{array}
$$

The estimates $\hat\epsilon_N$ and $\hat h_N$ are connected with the estimate $b_N$ of the parameter $b^*$ via this system of deterministic
algebraic equations. Therefore the rate of convergence $\hat\epsilon_N \to \epsilon$ and $\hat h_N \to h$ is determined by the rate of
convergence of $b_N$ to $b^*$ (which is exponential w.r.t. $N$). So we conclude that $\hat\epsilon_N \to\epsilon$ and $\hat h_N \to h$ $\mathbf{P}_{\epsilon}$-a.s.
as $N\to\infty$.

Theorem 2 is proved.

\end{document}